\documentclass[11pt]{amsart}

\usepackage{mathrsfs}
\usepackage{enumerate}
\newtheorem{theorem}{Theorem}[section]
\newtheorem{lemma}[theorem]{Lemma}

\theoremstyle{definition}
\newtheorem{definition}[theorem]{Definition}

\newtheorem{remark}{Remark}

\newcommand{\ra}{{\rightarrow}}
\newcommand{\lra}{{\longrightarrow}}
\newcommand{\eproof}{\hfill\rule{2.2mm}{3.0mm}}

\newcommand{\Proof}{\noindent {\bf Proof.~~}}
\newcommand{\D}{{\mathcal D}}

\newcommand{\R}{{\mathbb R}}
\newcommand{\Z}{{\mathbb Z}}

\newcommand{\M}{{\mathbb M}}
\newcommand{\ep}{\varepsilon}

\renewcommand{\eqref}[1]{(\ref{#1})}

\newcommand{\FT}[1]{\widehat{#1}}
\newcommand{\innerp}[1]{\langle {#1} \rangle}

\newcommand{\supp}{{\rm supp\,}}
\newcommand{\argmax}{{\rm argmax\,}}

\newcommand{\E}{{\mathcal E}}
\newcommand{\J}{{\mathcal J}}

\newcommand{\U}{\mathcal U}

\begin{document}

\title[The Regularity of Refinable Functions]
{{\Large\bf Letter to the Editor} \vspace{5mm} \\The Regularity of Refinable Functions}

\author{Yang Wang}
\thanks{Yang Wang was supported in part by the
       National Science Foundation grant DMS-0813750, DMS-08135022 and DMS-1043032.
       Zhiqiang Xu was supported  by NSFC grant 10871196 and by the Funds for Creative
       Research Groups of China (Grant No. 11021101).}
\address{Department of Mathematics  \\ Michigan State University\\
East Lansing, MI 48824, USA}
\email{ywang@math.msu.edu}

\author{Zhiqiang Xu}
\address{LSEC, Inst.~Comp.~Math., Academy of
Mathematics and System Science,  Chinese Academy of Sciences, Beijing, 100091, China}
\email{xuzq@lsec.cc.ac.cn}


\begin{abstract}
The regularity of refinable functions has been studied extensively in the past.
A classical result by Daubechies and Lagarias \cite{ingrid1} states that a compactly supported
refinable function in $\R$ of finite mask with integer dilation and translations
cannot be in $C^\infty$. A bound on the regularity based on the eigenvalues of certain
matrices associated with the refinement equation is also given. Surprisingly this fundamental classical result has not been proved in
the more general settings, such as in higher dimensions or when the dilation
is not an integer. In this paper we extend this classical result to the most general setting
for arbitrary dimension, dilation and translations.
\end{abstract}

\maketitle

\section{Introduction}
\setcounter{equation}{0}

 A {\em refinement equation} is a functional equation of the form
\begin{equation}\label{1.1}
f(x)= \sum_{d\in\D} c_d f(A x-d) %\mhsp  \sum_{j=0}^N c_j = |\det(A)|.
\end{equation}
where $\D\subset\R^n$ is a finite set, $c_d\neq 0$ for any $d\in\D$ and $A\in
\M_n(\R)$ is an $n\times n$ expanding matrix, i.e. all eigenvalues of
$A$ have $|\lambda|>1$. Since there are only
finitely many nonzero coefficients $c_d$, (\ref{1.1}) is often referred to as a refinement
equation with a {\em finite mask}. Here we shall refer to
a nontrivial function $f$ satisfying (\ref{1.1})
a {\em refinable function with dilation matrix
$A$ and translations $\D$}. In this paper, as in the vast majority
of studies in the literature, the focus is on
compactly supported refinable functions.

Refinement equations with finite masks
play a fundamental role in many applications such as
the construction of compactly
supported wavelets and in the study of subdivision schemes in CAGD. The regularity
of refinable functions is of great significance in those studies both in theory and in
applications. It has been studied extensively, including the seminal work by
Daubechies \cite{ingrid1989} which constructs compactly
supported refinable functions with orthogonal
integer translates of arbitrary regularity, leading to the fundamental class of
Daubechies wavelets. A more general study
by Daubechies and Lagarias \cite{ingrid1} establishes a classical result on the regularity of
a compactly supported refinable function in $\R$ of finite mask with integer dilation
and translations. It states that such a function cannot be in $C^\infty$, and it gives
bound on the regularity
based on the eigenvalues of certain matrices from the mask. The results in \cite{ingrid1}
have later been extended by several authors to obtain more refined regularity
estimations. In addition, using the same matrix eigenvalue technique,
one can extend the Daubechies-Lagarias result to refinable functions in $\R^n$ where
$A\in \M_n(\Z)$ and $\D \subset\Z^n$ (see Cabrelli, Heil and Molter \cite{CHM}).

Quite surprisingly, there have been very few results in terms of extending the classic
Daubechies-Lagarias result to the more general settings. This is perhaps due to
the fact that the matrix technique that has been effective for the integral case can
no longer be applied. Using techniques from number theory and harmonic analysis,
Dubickas and Xu \cite{DX}
prove that a refinable function in $\R$ with an arbitrary dilation $\lambda$ and
{\em integer} translations
cannot be in $C^\infty$. This appears to be the only generalization in this direction.
There have been other studies on the regularity of refinable functions with
non-integral dilations, e.g. in Dai, Feng and Wang \cite{DFW} on the decay rate of the
Fourier transform of a compactly supported
refinable function with arbitrary dilation and translations, and
in \cite{DFW2} by the same authors on refinable splines. There is also an extensive
literature on the absolute continuity of self-similar measures, which are
somewhat related to the study of regularity of refinable functions.
Nevertheless none of these studies directly address the extension of the
Daubechies-Lagarias result.

Many researchers in the community may have {\em assumed} that
the Daubechies-Lagarias result is valid in the general setting while in
reality other than those aforementioned special cases it has never been proved.
The general result turns out to be rather nontrivial to be established.
Our goal in this paper is to provide a short proof, thus establishing this important
classical result under the most general settings.
Our main theorem is:

\vspace{2mm}
\begin{theorem}  \label{theo-1.1}
     Let $f$ be a compactly supported refinable function of finite mask in $\R^n$.
     Then $f$ is not in $C^\infty(\R^n)$.
\end{theorem}

\vspace{2mm}

We shall prove the theorem in Section 2. In Section 3 we establish some upper bounds on
the regularity of compactly supported refinable functions.

\medskip

\noindent
{\bf Acknowledgement.} We thank Bin Han and Qiang Wu for helpful discussions.

\section{Proof of Main Theorem}
\setcounter{equation}{0}

Here in this section we prove  our main theorem.
We shall first investigate the support of $f$ by examining the attractor
of an iterated functions system (IFS) associated with a refinement equation and its
convex hull.
The IFS $\{\phi_d(x)=A^{-1}(x+d): ~d\in\D\}$ is referred to as the {\em IFS associated with
the refinement equation (\ref{1.1})}. By a well-known result of Hutchinson \cite{Hut81}
there is a unique compact set $T$ satisfying $T =\bigcup_{d\in\D} \phi_d(T)$.
The set $T$ is called the  {\em attractor} of the IFS $\{\phi_d: ~d\in \D\}$.
Let $\Phi(S):=\bigcup_{d\in\D} \phi_d(S)$ for any compact $S\subset\R^n$.
Then $T=\lim_{k\ra\infty} \Phi^k(S_0)$ in the Hausdorff
metric for any nonempty compact $S_0$. We shall let $\Omega:=\supp(f)$
denote the support of $f$.
It follows from the refinement equation (\ref{1.1}) that
$\Omega \subseteq \Phi(\Omega)$. By iterating it we obtain $\Omega\subset T$, where
$T$ is the attractor of the IFS $\{\phi_d: ~d\in\D\}$.

A key part of our tools involve the investigation of the convex hulls of various sets.
A point $z^*$ in a compact set $S$ is called
the {\em extremal point} of $S$ if
there is a unit vector $u\in\R^n$ such that
$z^*$ is the unique maximizer of $\innerp{u,x}$ for $x\in S$. In this case we shall
call $z^*$ the {\em extremal point of $S$ for the vector $u$.} These extremal points form
the extremal points of the convex hull of $S$. Let $\D_e\subseteq\D$ be the set of extremal points of the convex hull of $\D$. Somewhat related to this paper
are that the set of extremal
points of $T$ has been explicitly characterized in Strichartz and Wang \cite{StWa},
and furthermore it is shown in Dai and Wang \cite{DaWa09}
to be identical to the set of extremal points of $\Omega$.
%This result will be needed in Section 3.
%
%\begin{theorem}[\cite{DaWa09}]  \label{theo-2.1}
%     $\conv(T) =\conv(\Omega)$.
%\end{theorem}
%
%\medskip

Before proceeding further
we  first introduce some notations. For any $m\geq 1$ we define the map
$\pi_m: \D^m \lra \R^n$ by
\begin{equation}  \label{2.2}
    \pi_m\bigl([d_0, \dots, d_{m-1}]\bigl) := \sum_{j=0}^{m-1} A^j d_j.
\end{equation}
We let $\D_m:=\pi_m (\D^m)$, which is
$$
     \D_m :=\Bigl\{\sum_{j=0}^{m-1} A^j d_j:~[d_0, \dots, d_{m-1}]\in \D^m\Bigr\}.
$$
In general $\pi_m$ is not one-to-one. If $0\in\D$ then $\D_m \subseteq \D_{m+1}$.
We shall frequently consider the extremal points of $\D_m$ in this paper, and to
this end it is useful to introduce the set $\U$ of unit vectors defined by
$$
   \U := \Bigl\{u\in\R^n:~\|u\|=1~\mbox{and}~\innerp{u,d}\neq \innerp{u, e}~
   \mbox{for any distinct}~ d, e \in \bigcup_{m=1}^\infty \D_m\Bigr\}.
$$
Note that $\bigcup_{m=1}^\infty \D_m$ is a countable set so $\U$ is the whole
unit sphere in $\R^n$ minus a measure zero subset.

Now iterating the refinement equation (\ref{1.1}) we obtain
\begin{equation}  \label{2.3}
    f(x) = \sum_{v\in\D^m} c_v f\bigl(A^mx-\pi_m(v)\bigr)
         = \sum_{d\in\D_m} \tilde c_d f\bigl(A^mx-d\bigr),
\end{equation}
where $c_v = \prod_{j=0}^{m-1}c_{d_j}$ for $v=[d_0, \dots, d_{m-1}]$ and
$\tilde c_d = \sum_{v\in\D^m, \pi_m(v)=d} c_v$. Note that the support of the term
$f\bigl(A^mx-d\bigr)$ is $A^{-m}(\Omega+d)$.

\vspace{2mm}

\begin{lemma}  \label{lem-2.2}
    For any $u\in\U$ there is a unique $[d_0, d_1, \dots, d_{m-1}]\in\D^m$ such that
    $$
        \pi_m([d_0, d_1, \dots, d_{m-1}]) = \argmax_{d\in\D_m} \innerp{u, d}.
    $$
    Denote $d_m^u:= \pi_m([d_0, d_1, \dots, d_{m-1}])$.
    Consequently $\tilde c_{d_m^u} =\prod_{j=0}^{m-1}c_{d_j}$.
\end{lemma}
\Proof
The proof uses the same argument as in \cite{StWa}. We have
$\sum_{j=0}^{m-1}A^jd_j = \argmax_{d\in\D_m} \innerp{u, d}$.
Therefore for each $j$ we must have $A^jd_j =\argmax_{d\in\D} \innerp{u, A^jd}$. Furthermore,
because $u\in\U$ such $d_j$ must be unique, proving the lemma.
Of course in this case $\tilde c_{d_m^u} =\prod_{j=0}^{m-1}c_{d_j}$.
\eproof

\vspace{2mm}

\begin{definition}
    Let $\E$ be a subset of $\R^n$. Let $r_0>0$ and $u$ be a unit vector in $\R^n$.
    We say $\E$ has {\em $(u,r_0)$-isolated extremal point} if there exists
    a $z_0\in\E$ such that
    $\innerp{u,z_0} > \innerp{u,z}+r_0$ for all $z\in\E, z\neq z_0$.
\end{definition}

\vspace{2mm}

Intuitively if $\E$ has an $(u,r_0)$-isolated extremal point then the maximal point of
its projection onto the direction of $u$ is more than $r_0$ separated from the
other points in the projection. Our objective is to examine the sets $\D_m$ and
investigate whether it has certain isolation property. A key result is the following:

\vspace{2mm}

\begin{lemma}  \label{lem-2.3}
    Assume that there is an $r_0>0$ such that there exist infinitely many $m>0$
    such that $\D_m$ has $(u_m, r_0)$-isolated extremal point where
    $u_m\in\U$. Then a nontrivial compactly supported refinable function $f$
    satisfying (\ref{1.1}) cannot be in $C^\infty$.
\end{lemma}
\Proof
    Iterating the refinement equation leads to the new refinement equation (\ref{2.3}).
Now substitute $A^{-m}x$ for $x$ we obtain
\begin{equation}  \label{2.4}
    f(A^{-m}x) = \sum_{v\in\D^m} c_v f\bigl(x-\pi_m(v)\bigr)
         = \sum_{d\in\D_m} \tilde c_d f\bigl(x-d\bigr).
\end{equation}
Suppose $\D_m$ has $(u_m, r_0)$-isolated extremal point.
Let $z_m\in\Omega$ such that
$\innerp{u_m, z_m} = \sup_{x\in\Omega} \innerp{u_m, x}$. Clearly $z_m\in\partial\Omega$.
By the assumption
that $\D_m$ has $(u_m, r_0)$-isolated extremal point and $u_m\in\U$ we know that
\begin{equation}\label{eq:dmum}
    \innerp{u_m, d_m^{u_m}} > \innerp{u_m, d}+r_0
\end{equation}
for all other $d \neq d_m^{u_m}$ in $\D_m$. It follows that
$B_{r_0}(z_m)+d_m^{u_m}$ intersects the support of $f(x-d_m^{u_m})$ but is disjoint
from the support of all other $f(x-d)$ with $d \neq d_m^{u_m}$ in $\D_m$.

Now by assumption we have infinitely many $z_m \in \Omega$. Let
$z^*$ be a cluster point. Clearly there are infinitely many $m>0$ such that for
$r_1= r_0/2$,
$B_{r_1}(z^*)+d_m^{u_m}$ intersects the support of $f(x-d_m^{u_m})$ but is disjoint
from the support of all other $f(x-d)$ where $d_m^{u_m}\neq d\in\D_m$.
For any such $m$, setting $x=z+d_m^{u_m}$
in (\ref{2.4}) for $z\in B_{r_1}(z^*)$ yields
\begin{equation}  \label{2.5}
     f\bigl(A^{-m}z+A^{-m}d_m^{u_m}\bigr) = \tilde c_{d_m^{u_m}}f(z).
\end{equation}
Write $d_m^{u_m} =d_0+Ad_1+\cdots+A^{m-1}d_{m-1}$ where each $d_j\in\D$.  The equation (\ref{eq:dmum}) implies that $d_j\in \D_e\subset \D$.
By Lemma \ref{lem-2.2} we have $\tilde c_{d_m^{u_m}}= \prod_{j=0}^{m-1} c_{d_j}$. Hence
$|\tilde c_{d_m^{u_m}}|\geq b^m$ where $b = \min\{|c_d|: d\in\D_e\}$.
Now fix $x^*\in B_{r_1}(z^*)$ such that $| f(x^*)|>0$. Let
$$
     y_m = A^{-m}z_m+A^{-m}d_m^{u_m}, ~x_m =A^{-m}x^*+A^{-m}d_m^{u_m}.
$$
Since $z_m$ is on the boundary of $\Omega$, by (\ref{2.5}) we have $f(z_m) =0$. In fact,
since $z_m+d_m^{u_m}$ is an extremal point in $\Omega+\D_m$ for the vector $u_m$,
$y_m$ must be an extremal point of $A^{-m}(\Omega+\D_m)$ for the vector $(A^T)^mu_m$.
Hence $y_m$ must be on the boundary of $A^{-m}(\Omega+\D_m)$.
Now $A^{-m}(\Omega+\D_m)\supseteq \Omega$. It follows that $y_m$ must be either on
the boundary of $\Omega$ or not in the support of $f$ at all. In either case
if $f\in C^K(\R^n)$ then all $k$-th order
derivatives of $f$ must vanish at $y_m$ whenever $k \leq K$.
We also have $|f(x_m)|=|f(x_m)-f(y_m)| \geq b^m |f(x^*)|$.

We can derive a contradiction. Let $\tau<1$ be the spectral radius of $A^{-1}$. Then
$\|x_m-y_m\| = \|A^{-m}(x^*-z_m)\| \leq \tau^{-m}r_1$. Assume that $f$
is in $C_0^\infty(\R^n)$. Then for any $N>0$ and any
$y^*\in\partial\Omega$ we must have
$|f(y)-f(y^*)| =o(\|y-y^*\|^N)$ uniformly. In particular
$|f(x_m)-f(y_m)| =o(\|x_m-y_m\|^N)$. However,
\begin{equation} \label{2.6}
     \frac{|f(x_m)-f(y_m)|}{\|x_m-y_m\|^N}  \geq \frac{b^m|f(x^*)|}{\tau^{mN}\,r_1^N}
\end{equation}
for infinitely many $m$. By taking
$N$ large enough so that $\tau^N <b$ the right hand side of (\ref{2.6}) does not tend
to 0, a contradiction. Thus $f$ cannot be in
$C_0^\infty(\R^n)$.
\eproof

\vspace{2mm}

\begin{remark}
    The above proof actually gives a bound on the smoothness of $f$. A very crude
    bound that can be derived easily from (\ref{2.6}) is that
    if $f\in C^K(\R^n)$ then $K < \log b/\log \tau$ where
    $\tau$ is the spectral radius of $A^{-1}$ and $b = \min_{d\in\D_e} |c_d|$. The $\log b$
    part can be improved. In fact, let $\J$ be an infinite subset of indices $m$ such that
    $\D_m$ has $(u_m, r_0)$-isolated extremal point for each $m\in\J$. Then the proof
    shows that
    \begin{equation}  \label{2.7}
        K < \frac{\limsup_{m\in\J}\log |c_{d_m^{u_m}}|}{m\log\tau}.
    \end{equation}
    In many cases it allows us to obtain sharper upper bounds for
    the regularity of $f$.
\end{remark}

\vspace{2mm}

A key ingredient in our proof of the main theorem is the Borel-Cantelli Lemma, which states
that if $E_k$ is a sequence of events in some probability space and suppose that
$\sum_{k} {\rm Pr}(E_k) <\infty$ then ${\rm Pr}(\limsup_{k\rightarrow \infty} E_k) = 0$. In other words,
the probability that infinitely many of them occur is 0. We use the Borel-Cantelli Lemma
to prove the following key result.

\begin{theorem}  \label{theo-2.4}
     Let $A$ be an $n\times n$ dilation matrix and
     let $\D$ be a finite set in $\R^n$. Then there exists an $r_0>0$ and unit vectors
     $u_m\in\U$ such that $\D_m$ has $(u_m, r_0)$-isolated extremal
     point for infinitely many $m>0$.
\end{theorem}
\Proof
For $u\in\U$ let $e_m^u$ denote the element in $\D_m$ that
gives the second highest value of $\innerp{u, d}$ for $d\in\D_m$, i.e.
$$
     e_m^u = \argmax_{d\in\D_M\setminus\{d_m^u\}} \innerp{u, d}.
$$
If there is an $r_0>0$ such that
we can find infinitely many $m$ such that $\innerp{u, d_m^u-e_m^u} > r_0$ for some
$u=u_m\in\U$ we are done. Set
$$
     g_m(u) = \innerp{u, d_m^u-e_m^u}.% \mhsp\mbox{and}  \mhsp a_m = \sup_{u\in\U} g_m(u).
$$
We have $g_m(u)>0$ for all $m,u\in\U$. Assume the lemma is false.
Then $\lim_m \sup_{u\in\U} g_m(u) = 0$.
We shall derive a contradiction.

Assume that $d_m^u = \pi_m([d_0, d_1, \dots, d_{m-1}]) = \sum_{j=0}^{m-1}A^jd_j$.
We have already argued that each $d_j \in \D$ is the unique element in $\D$
satisfying $A^jd_j = \argmax_{d\in\D} \innerp{u, A^jd}$. Now assume that
$e_m^u = \sum_{j=0}^{m-1}A^je_j$. We claim that $d_j \neq e_j$ for one and only one
$0\leq j<m$. Assume $d_k \neq e_k$ and $d_l \neq e_l$
where $l\neq k$. Then the element $\tilde e_m^u =\sum_{j=0}^{m-1} A^j \tilde e_j$
where $\tilde e_j = e_j$ for all $j \neq l$ and $\tilde e_l = d_l$ will have the property
$$
      \innerp{u, d_m^u} > \innerp{u, \tilde e_m^u} > \innerp{u, e_m^u},
$$
contradicting the assumption that $\innerp{u, e_m^u}$ is the second largest. Thus there exists
a unique $0 \leq k <m$ and $d_k \neq e_k$ in $\D$ such that $d_m^u-e_m^u = A^k(d_k-e_k)$.
It follows that
$$
    g_m(u) =\innerp{u, A^k(d_k-e_k)}. % = \lambda^k \innerp{u, C(k\theta)(d_k-e_k)}.
$$
We shall denote $p_m(u):=k$
and $v_m(u):= d_k-e_k$. Thus we can rewrite the above equation as
$$
     g_m(u) = \innerp{u, A^{p_m(u)}v_m(u)}.
$$
Note that $\lim_m \sup_{u\in\U}g_m(u) =0$. So there exists a $C>0$ such that
$g_m(u) \leq C$ for all $m$ and $u\in\U$. Thus we have
\begin{equation}   \label{2.7}
    \innerp{u, A^{p_m(u)}v_m(u)} \leq C.
\end{equation}

Denote $\E := (\D-\D) \setminus\{0\}$, which is a finite set. Observe that $v_m(u) \in \E$.  For any $k>0$ define the set $E_k \subseteq \U$ by
$$
    E_k = \Bigl\{u\in\U:~\min_{v\in\E}\innerp{u, A^k v} \leq C \Bigr\}.
$$
Obviously, by (\ref{2.7}) if we pick $k=p_m(u)$ then $u\in E_k$. In particular if
$p_m(u)$ can take on infinitely many values then $u$ will be in infinitely many $E_k$.
We show that for almost all $u\in\U$, $p_m(u)$ can only take on finitely many values. To
see this we claim:

\noindent
{\bf Claim:}~~{\em There exists a constant $M>0$ such that
$\mu(E_k) \leq M \,\tau^k$, where $\mu$ is the normalized Hausdorff measure on the
unit sphere in $\R^n$ and $\tau$ is the spectral radius of $A^{-1}$.}

To prove the claim,
note that for each fixed unit vector $v_0$ the set
$F=\{u\in\U:~|\innerp{u, v_0}| \leq \ep\}$ has measure $\mu(F) \leq M_0\, \ep$
for some constant $M_0$ depending only on the dimension $n$. Now
for fixed $k$ and any $v\in\E$ set $R_k(v) =\|A^k v\|$ and $w_k(v) = A^k v/R_k(v)$.
It follows that the set
$$
    \{u\in\U:~|\innerp{u, w_k(v)}| \leq \ep\}
$$
has measure bounded by $M_0\,\ep$. Thus the set
$F_{k,v} =\{u\in\U:~|\innerp{u, A^k v}| \leq C\}$ has
$\mu(F_{k,v}) \leq M_0 C /R_m(v)$. Clearly,
$\min_{v\in\E} 1/R_m(v) \leq M_1\tau^m$ for some constant $M_1>0$.
Since $E_k$ is the union of $F_{k,v}$ with
where $v$ runs through $\E$, it follows that
$$
     \mu(E_k) \leq M_0\,C\,M_1\,L\,\tau^k
$$
where $L$ is the cardinality of $\E$. The claim follows
by setting $M = M_0\,C\,M_1\,L$.

Since $0<\tau<1$ we have $\sum_{k} \mu(E_k) <\infty$.
By the Borel-Cantelli Lemma, almost all $u\in\U$
belong to only finitely many $E_k$. In other words, for almost all $u\in \U$ the set
$\{p_m(u)\}$ is a finite set. Thus taking any such $u\in \U$, the sequence
$g_m(u) = \innerp{u, A^{p_m(u)}v_m(u)}$
can take on only finitely many values. Furthermore we already know that
$g_m(u) \neq 0$. This contradicts the assumption that $\lim_m g_m(u) =0$.
The lemma is now proved.
\eproof

\vspace{2mm}

Theorem \ref{theo-1.1} now follows easily from
Lemma \ref{lem-2.3} and Theorem \ref{theo-2.4}.

\section{Regularity Upper Bounds}
\setcounter{equation}{0}

In proving our main theorem in the previous section we have in fact already established
upper bounds for the regularity of refinable functions, or at least ways to estimate
such bounds. We have already remarked in Section 2 that a very simple but crude
bound for a compactly supported refinable function satisfying (\ref{1.1}) is
\begin{equation}  \label{3.0}
        K < \frac{\log b}{\log\tau},
\end{equation}
where $b=\min\{|c_d|: d\in\D_e\}$ and $\tau$ is the spectral radius of $A^{-1}$. We
establish some more refined bounds here.
\smallskip

\begin{theorem}  \label{theo-3.2}
     Let $f$ be a compactly supported refinable function in $\R^n$ satisfying the refinement
     equation (\ref{1.1}). Let $w\in\R^n$ be a unit eigenvector of $A^T$ corresponding to
     a real eigenvalue $\lambda$ of $A^T$. Assume that $\D$ has an
     extremal point $d_w\in\D$ for the vector $w$ and $f\in C^K(\R^n)$, $K \geq 0$. Then
     \begin{equation}  \label{3.1}
         K < \frac{\log |c_{d_w}|}{\log|\lambda|^{-1}}.
    \end{equation}
\end{theorem}
\Proof Since $d_w$ is extremal in $\D$ for the vector $w$ we can find an $r_0>0$
such that $r_0< \innerp{w, d_w} - \max_{\D\setminus\{d_w\}}\innerp{w,d}$. With $w$ being
an eigenvector of $A^T$ it is
straightforward to see that $d_m^w=d_w+Ad_w+ \cdots+A^{m-1}d_w$ is the
extremal point of $\D_m$ for the vector $w$, and it gives $\D_m$ a $(w, r_0)$-isolated
extremal point for each $m$. Let $z^* = \argmax_{x\in\Omega} \innerp{w,x}$. Then for
any $d\neq d_m^w$ in $\D_m$ and $x\in\Omega$ we must have
\begin{equation} \label{3.2}
     \innerp{w,z^*+d_m^w} > \innerp{w, x+d} + r_0.
\end{equation}
Thus if $z\in\Omega$ such that $\innerp{w, z^*-z} \leq r_0$ then $y=d_m^w+z$ is
not in $d+\Omega$ for all $d\neq d_m^w$ in $\D_m$, and hence it satisfies $f(y-d)=0$.

Note that $f$ is not identically
0 in any neighborhood of $z^*$ since $z^*$ is in the support of $f$. Pick an
$x^*\in B_{r_0}(z^*)$ such that $f(x^*) \neq 0$. Now let $u\in\R^n$ be a unit
$\lambda$-eigenvector of $A$ such that $\innerp{w,u}\geq 0$. Consider the
set $E \subset \R^+$,
$$
    E = \bigl\{t \geq 0: ~x^*+ t u \in\Omega\bigr\}.
$$
Let $t_0 = \max E$ and $z_0 = x^*+t_0u$. Clearly $z_0\in\partial\Omega$.
Furthermore $z_0$ satisfying (\ref{3.1}) since
$$
   \innerp{w, z^*-z_0} = \innerp{w,z^*-x^*} - t\innerp{w,u}
      \leq \innerp{w,z^*-x^*} \leq r_0.
$$
Hence $d_m^w+z_0$ is not in $d+\Omega$ for all $d\neq d_m^w$ in $\D_m$. But
$z_0+ d_m^w\in\partial(\Omega+d_m^w)$. It follows that
$z_0+ d_m^w \in \partial(\Omega+\D_m)$. Set $y_m = A^{-m}(z_0+d_m^w)$,
which must be on the boundary of  $A^{-m}(\Omega+\D_m)$.
Now $A^{-m}(\Omega+\D_m)\supseteq \Omega$, which implies that $y_m$ must be either on
the boundary of $\Omega$ or not in the support of $f$ at all. In either case
since $f\in C^K(\R^n)$ all $k$-th order
derivatives of $f$ must vanish at $y_m$ whenever $k \leq K$.

We can now bound the regularity $K$ of $f$.
Denote $x_m = A^{-m}(x^*+d_m^w)$. We have
$f(x_m) = \tilde c_{d_m^w} = c_{d_w}^m f(x^*)$, and hence
$|f(x_m)-f(y_m)|=|f(x_m)| \geq |c_{d_w}|^m |f(x^*)|$. Observe that
$\|x_m-y_m\| = \|A^{-m}(x^*-z_0)\| =|\lambda|^{-m}\|x^*-z_0\|$ since
$x^*-z_0 = t_0 u$ is a $\lambda$-eigenvector of $A$.
Since $f\in C_0^K(\R^n)$, for any $y^*\in\partial\Omega$ and $y\in \Omega$ we must have
$|f(y)-f(y^*)| =o(\|y-y^*\|^K)$ uniformly. In particular
$|f(x_m)-f(y_m)| =o(\|x_m-y_m\|^K)$. However,
\begin{equation} \label{3.3}
     \frac{|f(x_m)-f(y_m)|}{\|x_m-y_m\|^K}
          = \frac{|c_{d_w}|^m |f(x^*)|}{\|x^*-z_0\|^K\,|\lambda|^{-mK}}.
\end{equation}
The right hand side of (\ref{3.3}) will tend
to 0 as $m \ra \infty$ only if $K < \frac{\log |c_{d_w}|}{\log|\lambda|^{-1}}$.
This proves the theorem.
\eproof

\smallskip
\begin{remark}
It is well known that without the compactly supported assumption $f$ can in fact be $C^\infty$. The simplest example is $f(x) = x$, which satisfies $f(x) = \frac{1}{2}f(2x)$. The study
of refinable functions has largely imposed the additional condition
$\sum_{d\in\D} c_d = |\det(A)|$, which stems from many applications such as wavelets.
Under this condition, a refinement equation with finite mask
has up to a scalar multiple a unique
compactly supported distribution solution. If we restrict to only solutions in
$L^1(\R)$ then the uniqueness result also
holds without the additional sum of coefficients condition in the one dimension, provided
such a solution exists \cite{ingrid1}. Surprisingly, just like the regularity
result before this paper, this result has not been
established for higher dimensions in the general setting.
\end{remark}

\smallskip

\begin{remark}
With the additional condition $\sum_{d\in\D} c_d = |\det(A)|$ it is well known that
a compactly supported refinable distribution $f$ must have $\FT f(0) \neq 0$. This fact
can be combined with the projection method in \cite{han1} to yield a slightly
less tedious proof of
Theorem \ref{theo-3.2}. Without the additional condition, however, one problem we
cannot overcome is to show that the projection is nontrivial.
\end{remark}

\bigskip

\bibliographystyle{amsplain}

\end{document}